\documentclass[11pt]{article}

\usepackage[margin=1in]{geometry}
\usepackage{amsmath,amssymb,amsthm}
\usepackage[ruled,linesnumbered,vlined]{algorithm2e}
\usepackage{enumitem}
\usepackage{booktabs}
\usepackage[hidelinks]{hyperref}

\newcommand{\IZ}{\mathbb{Z}}

\newcommand{\IR}{\mathbb{R}}

\newcommand{\abs}[1]{\left\lvert #1\right\rvert}
\newcommand{\norm}[1]{\left\lVert #1\right\rVert}
\newcommand{\fr}[1]{\left\{#1\right\}}          
\newcommand{\seq}[1]{\left\langle #1\right\rangle}

\newtheorem{theorem}{Theorem}[section]
\newtheorem{proposition}[theorem]{Proposition}
\newtheorem{lemma}[theorem]{Lemma}
\newtheorem{corollary}[theorem]{Corollary}
\theoremstyle{definition}
\newtheorem{definition}[theorem]{Definition}
\theoremstyle{remark}
\newtheorem{remark}[theorem]{Remark}

\numberwithin{equation}{section}
\SetKwInOut{Input}{input}\SetKwInOut{Output}{output}

\title{Ostrowski Numeration, Alternating Expansions\\ and Inhomogeneous Linear Approximation}
\author{Avraham Bourla}
\date{}

\begin{document}
\maketitle

\begin{abstract}
\noindent Given an irrational slope $\alpha$ and a real intercept $\beta$, the inhomogeneous linear approximation problem asks
for integers $m$ for which $\norm{m\alpha-\beta}$ is small compared with $\abs{m}$. Building on Ostrowski's numeration system,
we give a self-contained treatment of its alternating (dual) counterpart, whose weights are the signed convergent
denominators $(-1)^kq_k$, and show that with $n$ digits it represents exactly the integers of an explicit interval of length
$q_n+q_{n-1}$, compared with $q_n$ for Ostrowski's system. We pair the absolute and the alternating integer systems with,
respectively, Ostrowski's classical alternating expansion and an absolute expansion of real numbers in the linear
approximation coefficients $\theta_k=q_k\alpha-p_k$. These expansions yield explicit algorithms that produce,
for the one--sided and the two--sided (forward and backward) variants of the problem, sequences of approximates $A_n$ with
$\abs{A_n}\le q_n$ and error at most $\abs{\theta_{n-1}}$. For the one--sided variants we prove optimality: the $n$--th term
is the best one--sided approximate within the $n$--th range of the alternating system.
\end{abstract}

\noindent\textbf{Keywords:} inhomogeneous Diophantine approximation, continued fractions, Ostrowski numeration,
alternating numeration systems.\\
\textbf{MSC 2020:} 11J70, 11A63, 11J20, 11K16.

\section{Introduction}\label{sec_intro}

Fix an irrational slope $\alpha\in(0,1)$ and an intercept $\beta\in[0,1)$ (there is no loss of generality, since only the
fractional parts matter). Writing $\norm{r}$ for the distance from $r\in\IR$ to the nearest integer, the
\textbf{inhomogeneous linear approximation problem} asks for integers $m$ that make the error $\norm{\beta-m\alpha}$ small
while keeping $\abs{m}$ small; we borrow the terms \emph{slope} and \emph{intercept} from \cite{R}. The homogeneous case
$\beta=0$ is completely solved by the continued fraction expansion of $\alpha$: the convergent denominators $q_k$ are
exactly the integers $m\ge1$ with $m\norm{m\alpha}<1$ (Theorem 16 in \cite{K}). For $\beta\ne0$ no such normality
condition can hold in general (Theorem 26 in \cite{K}), and a quantitative relaxation is needed.

The following example illustrates why the inhomogeneous problem requires its own tools. Take $\alpha:=\pi-3$ and
$\beta:=e-2$. As homogeneous approximates, $2$ and $-2$ both have error $\norm{2\pi}\approx0.28$, and since this is
close to $1-\fr{e}$, the integer $-2$ is an excellent backward $(\alpha,\beta)$--approximate. A direct computation shows that
$\norm{e+2\pi}<\norm{e-n\pi}$ for every $0\le n\le22251$: one has to go as far as $n=22252$ to find a forward approximate
that is more accurate than $-2$.

\paragraph{Prior work.} Expansions of real numbers with respect to the sequence $\seq{\theta_k}$ go back to
Ostrowski \cite{Ost}, who also introduced the numeration system on the convergent denominators that bears his name; see
\cite{RS} and the surveys \cite{Ber,BI} for its many applications. Inhomogeneous approximation is classical
(see e.g.\ \cite[Chapter~III]{Cas} and \cite{Bug}), and algorithms describing best inhomogeneous approximations by means of
Ostrowski--type expansions of the intercept were given in \cite{NST,Kom}; see also \cite{Ito,IN} for the associated
dynamics. The alternating integer system of Section~\ref{sec_numeration} is the dual Ostrowski system of
\cite[Proposition~6.4.19]{PF}; for $a_k\equiv1$ it is Bunder's representation of the integers by negative Fibonacci numbers
\cite{Bun}, and for constant partial quotients it is studied in \cite{Fok}. Cabanillas \cite{Cab} introduced a variant of
Ostrowski's system in which an integer $n$ and the fractional part $\fr{n\alpha}$ share the same digit string, coding all
of $\IZ$ and $[0,1)$, and used it to describe the best one--sided approximations of an intercept by truncating its digits
\cite[\S4.3]{Cab}. The alternating expansion of Theorem~\ref{gamma} is Ostrowski's classical expansion of the intercept,
see \cite{Ost,RS,Ber}. The present paper develops the subject from the point of view of numeration systems.

\paragraph{Contributions.} The main results are the following.
\begin{enumerate}[label=(\arabic*)]
\item A self-contained, constructive treatment of the \emph{alternating integer system} with the signed weights
$q_k^*:=(-1)^kq_k$ (Theorem~\ref{thm_T}). Besides existence and uniqueness, which are known \cite{PF}, we show that the
strings of length at most $n$ cover exactly an explicit interval $I_n^*$ of $q_n+q_{n-1}$ integers, and that the length of
the representation of $T$ is the unique $n$ with $T\in I_n^*\setminus I_{n-1}^*$.
\item A \emph{duality} between integer systems and real expansions: the alternating integer system pairs with an
\emph{absolute} expansion $\beta=\sum b_k\abs{\theta_{k-1}}$ (Theorem~\ref{beta}), and Ostrowski's absolute system pairs
with the classical \emph{alternating} expansion $\gamma=\sum c_k\theta_{k-1}$ (Theorem~\ref{gamma}), for which we give a
self-contained proof by tail intervals and an iterative algorithm with the correct endpoint convention.
\item \emph{Explicit approximation sequences} for the four signed variants of the problem, obtained by truncating these
expansions (Corollaries~\ref{total_under}, \ref{total_over}, \ref{positive_two_sided} and \ref{negative_two_sided}), and an
\emph{optimality theorem} for the one--sided variants within the ranges $I_n^*$, which contain negative integers
(Theorem~\ref{T}); compare \cite[\S4.3]{Cab} for best one--sided approximations by non--negative integers.
\end{enumerate}
All constructions are algorithmic (Algorithms~\ref{algo_S}--\ref{algo_gamma}).

\paragraph{Organisation.} Section~\ref{sec_prelim} collects the facts on continued fractions and on the sequence
$\seq{\theta_k}$ that we use. Section~\ref{sec_problem} states the homogeneous and inhomogeneous problems and the notion of a
general solution. Section~\ref{sec_numeration} treats the two integer systems and Section~\ref{Irrational Base Expansions} the
two real expansions and the approximation theorems. Section~\ref{sec_remarks} discusses rational slopes, where the same
identities give modular multiplication with the quotient, and open questions. The long proof of Theorem~\ref{thm_T} is in
Appendix~\ref{app_T}.

\paragraph{Notation.} $\lfloor r\rfloor$ and $\lceil r\rceil$ denote floor and ceiling, $\fr{r}:=r-\lfloor r\rfloor$ the
fractional part and $\norm{r}:=\min\{\fr{r},1-\fr{r}\}$. Digit strings are written $\seq{d_k}_1^n$; an extended index
$\ell\in\IZ^*_{\ge0}:=\IZ_{\ge0}\cup\{\infty\}$ allows infinite strings, with $\infty\pm1:=\infty$. The parity of $n$ is
$\rho_n:=\lceil n/2\rceil-\lfloor n/2\rfloor\in\{0,1\}$.

\section{Continued fractions and linear approximation coefficients}\label{sec_prelim}

\subsection{Numeration systems}
A \textbf{numeration system} is a bijection between $\IZ_{\ge0}$ and a set of finite strings of non--negative digits
(the \textbf{admissible} strings); the value of a string is its dot product with a fixed \textbf{weight vector}, as
$\seq{10^k}$ for the decimal system. An \textbf{integer system} is the analogous bijection with $\IZ$. Negating the ratio
of the weights gives the base $-10$ integer system, whose $n$--digit strings cover the interval
$[-9090,909]\cap\IZ$ for $n=4$: the same $10^4$ integers as the absolute system, merely shifted. Section~\ref{sec_numeration}
shows that for continued fraction weights the alternating system is genuinely larger than the absolute one.

\subsection{The continued fraction expansion}
Given $\alpha\in[0,1)$ we expand it with the Gauss map:

\begin{algorithm}[H]
\Input{the initial seed $\alpha\in[0,1)$}
\Output{the digit string $\seq{a_k}_1^\ell\subseteq\IZ_{\ge1}$}
set $\alpha_0:=\alpha;\ \ell:=\infty;\ k:=1$\;
\While{$\alpha_{k-1}>0$}{
  set $a_k:=\lfloor1/\alpha_{k-1}\rfloor$\;
  set $\alpha_k:=\fr{1/\alpha_{k-1}}=1/\alpha_{k-1}-a_k\in[0,1)$\;
  set $k:=k+1$\;
}
set $\ell:=k-1$\;
\If{$\ell\ge2$ and $a_\ell=1$}{set $\ell:=\ell-1$;\ $a_\ell:=a_\ell+1$\;}
\caption{the continued fraction expansion}\label{algo_cf}
\end{algorithm}

\noindent From line~3 we have
\begin{equation}\label{a_k}
a_k\alpha_{k-1}\le1<(a_k+1)\alpha_{k-1},\qquad 1\le k<\ell+1,
\end{equation}
and line~4 gives $\alpha_{k-1}=(a_k+\alpha_k)^{-1}$, hence the continued fraction
$\alpha=[0;a_1,a_2,\dots]$. Its $n$--th \textbf{convergent} $p_n/q_n:=[0;a_1,\dots,a_n]$ satisfies the classical
facts (see \cite{K})
\begin{equation}\label{convergents_bound}
\abs{\alpha-\frac{p_k}{q_k}}<\frac1{q_k^2},\qquad
\frac{p_{2k}}{q_{2k}}<\alpha<\frac{p_{2k+1}}{q_{2k+1}}\quad(k\ge0),
\end{equation}
\begin{equation}\label{p_k_recursion}
p_{-1}:=1,\quad p_0:=0,\quad p_k=a_kp_{k-1}+p_{k-2},\qquad k\ge1,
\end{equation}
\begin{equation}\label{q_k_recursion}
q_{-1}:=0,\quad q_0:=1,\quad q_k=a_kq_{k-1}+q_{k-2},\qquad k\ge1.
\end{equation}
The expansion is finite exactly when $\alpha$ is rational and is unique once expansions ending with the digit $1$ are
excluded (line~8). \emph{From now on $\alpha\in(0,1)$ is irrational}, except in Section~\ref{sec_remarks}.

Iterating \eqref{q_k_recursion} and using $q_0=1$ gives
\begin{equation}\label{q_n_even_odd}
q_n=\begin{cases}1+\sum_{k=1}^{n/2}a_{2k}q_{2k-1}, & n\text{ even},\\[2pt] \sum_{k=0}^{(n-1)/2}a_{2k+1}q_{2k}, & n\text{ odd}.\end{cases}
\end{equation}
We write
\begin{equation}\label{p_k^*_q_k^*}
p_k^*:=(-1)^kp_k,\qquad q_k^*:=(-1)^kq_k,\qquad k\ge-1,
\end{equation}
so that $p_k^*=p_{k-2}^*-a_kp_{k-1}^*$ and $q_k^*=q_{k-2}^*-a_kq_{k-1}^*$ for $k\ge1$. When $a_k=1$ for all $k$,
$\alpha=\phi:=(\sqrt5-1)/2$ and $\seq{q_k}$ is the Fibonacci sequence; when $a_k=2$ for all $k$, $\alpha=\sqrt2-1$ (the
\emph{silver section}) and $\seq{q_k}_0^4=\seq{1,2,5,12,29}$.

\subsection{The sequence of linear approximation coefficients}
For the convergents of $\alpha$ define
\begin{equation}\label{theta_def}
\theta_k:=q_k\alpha-p_k,\qquad k\ge-1.
\end{equation}
By \eqref{p_k_recursion} and \eqref{q_k_recursion},
\begin{equation}\label{theta_k_recursion}
\theta_{-1}=-1,\qquad\theta_0=\alpha,\qquad\theta_k=\theta_{k-2}+a_k\theta_{k-1},\qquad k\ge1,
\end{equation}
and an induction with line~4 of Algorithm~\ref{algo_cf} gives
\begin{equation}\label{theta_induction}
\theta_k=-\theta_{k-1}\alpha_k=(-1)^k\alpha_0\alpha_1\cdots\alpha_k,\qquad k\ge0.
\end{equation}
Hence $\seq{\abs{\theta_k}}$ decreases strictly to $0$ and alternates in sign,
\begin{equation}\label{theta_alt}
\abs{\theta_k}=(-1)^k\theta_k,\qquad\abs{\theta_k}=q_k^*\alpha-p_k^*,\qquad k\ge-1,
\end{equation}
\begin{equation}\label{abs_theta_k_recursion}
\abs{\theta_{-1}}=1,\qquad\abs{\theta_0}=\alpha,\qquad\abs{\theta_k}=\abs{\theta_{k-2}}-a_k\abs{\theta_{k-1}},\qquad k\ge1.
\end{equation}
Since $2\abs{\theta_1}=2\alpha_1/(a_1+\alpha_1)<1$, we have
\begin{equation}\label{norm_abs}
\norm{\theta_k}=\abs{\theta_k},\qquad k\ge1.
\end{equation}
From \eqref{a_k} and \eqref{theta_induction}, $a_k\abs{\theta_{k-1}}/\abs{\theta_{k-2}}=a_k\alpha_{k-1}$, so
\begin{equation}\label{theta_a_k}
a_k\abs{\theta_{k-1}}<\abs{\theta_{k-2}}<(a_k+1)\abs{\theta_{k-1}},\qquad k\ge1.
\end{equation}

\begin{proposition}\label{conv}
The series $\sum_{k\ge0}\abs{\theta_k}$ converges. Consequently $\sum_{k\ge1}a_k\abs{\theta_{k-1}}$ and
$\sum_{k\ge1}a_k\theta_{k-1}$ converge absolutely.
\end{proposition}
\begin{proof}
By line~4 of Algorithm~\ref{algo_cf}, $\alpha_k=1/(a_{k+1}+\alpha_{k+1})$, hence
$\alpha_k\alpha_{k+1}=\alpha_{k+1}/(a_{k+1}+\alpha_{k+1})\le\alpha_{k+1}/(1+\alpha_{k+1})<\tfrac12$. By \eqref{theta_induction},
$\abs{\theta_{k+1}}=\abs{\theta_{k-1}}\alpha_k\alpha_{k+1}<\tfrac12\abs{\theta_{k-1}}$ for $k\ge0$, so the series converges
geometrically. The second statement follows from $a_k\abs{\theta_{k-1}}<\abs{\theta_{k-2}}$, see \eqref{theta_a_k}.
\end{proof}

\noindent Telescoping \eqref{abs_theta_k_recursion} and \eqref{theta_k_recursion} gives, for $n\ge1$,
\begin{equation}\label{partial_abs_theta}
\sum_{k=1}^na_k\abs{\theta_{k-1}}=1+\alpha-\abs{\theta_{n-1}}-\abs{\theta_n},\qquad
\sum_{k=1}^na_k\theta_{k-1}=1-\alpha+\theta_{n-1}+\theta_n,
\end{equation}
and letting $n\to\infty$,
\begin{equation}\label{abs_rep}
1+\alpha=\sum_{k=1}^\infty a_k\abs{\theta_{k-1}},\qquad
\sum_{k=n+1}^\infty a_k\abs{\theta_{k-1}}=\abs{\theta_{n-1}}+\abs{\theta_n},
\end{equation}
\begin{equation}\label{cond_rep}
1-\alpha=\sum_{k=1}^\infty a_k\theta_{k-1}.
\end{equation}
Taking half the difference and half the sum of \eqref{abs_rep} and \eqref{cond_rep},
\begin{equation}\label{self_rep}
\alpha=\sum_{k=1}^\infty a_{2k}\abs{\theta_{2k-1}}=-\sum_{k=1}^\infty a_{2k}\theta_{2k-1},\qquad
1=\sum_{k=0}^\infty a_{2k+1}\theta_{2k}=\sum_{k=0}^\infty a_{2k+1}\abs{\theta_{2k}}.
\end{equation}
Finally, iterating the recursions from the index $n-1$ gives the tail expansions
\begin{equation}\label{abs_theta_tail_odd_even}
\abs{\theta_{n-1}}=\sum_{j\ge0}a_{n+1+2j}\abs{\theta_{n+2j}},\qquad
-\theta_{n-1}=\sum_{j\ge0}a_{n+1+2j}\theta_{n+2j},\qquad n\ge0.
\end{equation}

\section{The approximation problems}\label{sec_problem}

\subsection{The solved homogeneous case}
An integer $m\ne0$ is an $\alpha$--approximate with \textbf{two--sided error} $\norm{m\alpha}$; it is \textbf{normal} if
\begin{equation}\label{normal}
\abs{m}\cdot\norm{m\alpha}<1.
\end{equation}
By \eqref{convergents_bound}, $q_k\norm{q_k\alpha}\le q_k\abs{q_k\alpha-p_k}<1$, so every $q_k$ is normal, and by Theorem~16
of \cite{K} there are no other normal positive approximates; since $\norm{-r}=\norm{r}$, the normal negative approximates are
the $-q_k$. Call $m$ \textbf{over} if $\fr{m\alpha}<1/2$ and \textbf{under} otherwise. By \eqref{convergents_bound}, $q_k$ is
over exactly when $k$ is even and $-q_k$ exactly when $k$ is odd, so the normal over and under approximates are the terms of
$\seq{q_k^*}$ and $\seq{-q_k^*}$ respectively. The quadratic coefficients $q_k\abs{\theta_k}$ have a rich structure of their
own, see \cite{B1,DK}.

\subsection{The inhomogeneous problem}
For the intercept $\beta\in(0,1)$ the integer $m$ (now including $m=0$) is an $(\alpha,\beta)$--approximate with error
$\norm{\beta-m\alpha}$. We call the non--negative, non--positive and arbitrary approximates \textbf{forward},
\textbf{backward} and \textbf{total}, and say that $m$ is \textbf{under} $\beta$ if $\fr{m\alpha}\in(\beta-\tfrac12,\beta]
\pmod 1$ and \textbf{over} $\beta$ otherwise.

\begin{definition}\label{def_solution}
Let $\ell\in\IZ^*_{\ge0}$. A sequence $\seq{A_n}_1^\ell$ of $(\alpha,\beta)$--approximates is a \textbf{general solution} if
\begin{enumerate}[label=(\roman*)]
\item $\abs{A_n}\le q_n$ for $1\le n<\ell+1$;
\item $0<\norm{\beta-A_n\alpha}\le\abs{\theta_{n-1}}$ for $1\le n<\ell$;
\item $\norm{\beta-A_\ell\alpha}=0$ if $\ell<\infty$, and $\norm{\beta-A_n\alpha}\to0$ if $\ell=\infty$.
\end{enumerate}
\end{definition}
\noindent Thus the $n$--th term is as large as the $n$--th homogeneous solution $q_n$ while its error is bounded by the error
$\abs{\theta_{n-1}}$ of the previous one. For the one--sided solutions of Corollaries~\ref{total_under} and~\ref{total_over}
the inequality in (ii) is strict; for the two--sided solutions it is strict whenever $\beta\notin\IZ+\alpha\IZ$
(Remark~\ref{rem_strict}).

\section{Irrational base numeration systems}\label{sec_numeration}

\subsection{Ostrowski's absolute system}
A digit string $\seq{c_k}_1^n$ is \textbf{left--admissible} if it satisfies the Markov conditions
\begin{equation}\label{left--admissible}
\text{(i) } 0\le c_1\le a_1-1;\qquad\text{(ii) } 0\le c_k\le a_k\ (2\le k\le n);\qquad
\text{(iii) } c_k=a_k\implies c_{k-1}=0\ (2\le k\le n).
\end{equation}

\begin{theorem}[Ostrowski \cite{Ost}]\label{thm_S}
Every $S\in\IZ_{\ge0}$ has a unique representation $S=\sum_{k=1}^nc_kq_{k-1}$ with $\seq{c_k}_1^n$ left--admissible and
$c_n\ge1$; its index satisfies $q_{n-1}\le S\le q_n-1$. In particular the left--admissible strings of length at most $n$
represent exactly the integers $0\le S\le q_n-1$.
\end{theorem}

\noindent The theorem is classical; we include the greedy algorithm and a proof for completeness.

\begin{algorithm}[H]
\Input{the counting number $S\in\IZ_{\ge0}$}
\Output{the left--admissible digit string $\seq{c_k}_1^n$}
set $S_0:=S;\ m:=0$\;
\While{$S_m\ge1$}{
  let $n_m$ be such that $q_{n_m-1}\le S_m\le q_{n_m}-1$\;
  set $c_{n_m}:=\lfloor S_m/q_{n_m-1}\rfloor$\;
  set $S_{m+1}:=S_m-c_{n_m}q_{n_m-1}$\;
  set $m:=m+1$\;
}
set $M:=m;\ n:=n_0;\ c_k:=0$ for all $k\notin\{n_m\}_{m=0}^{M-1}$\;
\caption{the base--$\alpha$ numeration system}\label{algo_S}
\end{algorithm}

\begin{proof}
Line~4 gives $c_{n_m}q_{n_m-1}\le S_m<(c_{n_m}+1)q_{n_m-1}$, hence $S_{m+1}<q_{n_m-1}$ and the indices strictly decrease,
$0\le n_{M}<\dots<n_1<n_0=n$; the iteration stops with $S_M=0$ and $S=\sum_kc_kq_{k-1}$. (If $a_1=1$ then $q_1=q_0=1$ and
for $S=1$ line~3 picks $n_m=2$.) Admissibility: if the last step has $n_{M-1}=1$, then $c_1=S_{M-1}\le q_1-1=a_1-1$, which is
(i). If $c_{n_m}\ge a_{n_m}+1$ then $S_m<q_{n_m}=a_{n_m}q_{n_m-1}+q_{n_m-2}<(a_{n_m}+1)q_{n_m-1}\le S_m$, a contradiction,
which gives (ii). If $c_k=a_k$ and $c_{k-1}\ge1$ for some $k\ge2$, the two digits were produced in consecutive steps
$m$ (index $k$) and $m+1$ (index $k-1$), and $S_m\ge a_kq_{k-1}+q_{k-2}=q_k$, contradicting $S_m\le q_{n_m}-1=q_k-1$; this
gives (iii). Uniqueness: suppose two left--admissible strings represent $S$; both have the same index $n$ by the bound on
$S$, and if $j$ is the largest index where they differ, say $c_j<\hat c_j$, then
$\sum_{k<j}c_kq_{k-1}\le q_{j-1}-1$ (the bound applied to the truncated string) yields
$\sum_{k\le j}c_kq_{k-1}<(c_j+1)q_{j-1}\le\sum_{k\le j}\hat c_kq_{k-1}$, so the two values differ. The bound
$\sum_{k=1}^nc_kq_{k-1}\le q_n-1$ for left--admissible strings follows by induction on $n$ from \eqref{q_k_recursion}
and the conditions (i)--(iii).
\end{proof}

\noindent For $\alpha=\phi$ this is Zeckendorf's theorem. For the silver section $\alpha=\sqrt2-1$ ($a_k=2$,
$\seq{q_k}_0^3=\seq{1,2,5,12}$) the strings of length $4$ represent $0,\dots,q_4-1=28$; Table~\ref{tab_abs} lists them.

\begin{table}[h]\centering\small
\begin{tabular}{r|cccc||r|cccc}
\toprule
$N$ & $c_4$ & $c_3$ & $c_2$ & $c_1$ & $N$ & $c_4$ & $c_3$ & $c_2$ & $c_1$\\
\midrule
1&0&0&0&1 & 15&1&0&1&1\\ 2&0&0&1&0 & 16&1&0&2&0\\ 3&0&0&1&1 & 17&1&1&0&0\\ 4&0&0&2&0 & 18&1&1&0&1\\
5&0&1&0&0 & 19&1&1&1&0\\ 6&0&1&0&1 & 20&1&1&1&1\\ 7&0&1&1&0 & 21&1&1&2&0\\ 8&0&1&1&1 & 22&1&2&0&0\\
9&0&1&2&0 & 23&1&2&0&1\\ 10&0&2&0&0 & 24&2&0&0&0\\ 11&0&2&0&1 & 25&2&0&0&1\\ 12&1&0&0&0 & 26&2&0&1&0\\
13&1&0&0&1 & 27&2&0&1&1\\ 14&1&0&1&0 & 28&2&0&2&0\\
\bottomrule
\end{tabular}
\caption{The base--$\alpha$ system for $\alpha=\sqrt2-1$, weights $q_3,q_2,q_1,q_0=12,5,2,1$.}\label{tab_abs}
\end{table}

\subsection{The alternating integer system}
A digit string $\seq{b_k}_1^n$ is \textbf{right--admissible} if
\begin{equation}\label{right--admissible}
\text{(i) } 0\le b_k\le a_k\ (1\le k\le n);\qquad\text{(ii) } b_k=a_k\implies b_{k+1}=0\ (1\le k\le n-1).
\end{equation}
For $n\ge0$ define the integer intervals
\begin{equation}\label{I_n^*}
I_n^*:=\bigl[\,1-q_{n-\rho_n},\ q_{n-\rho_{n-1}}\bigr]\cap\IZ,
\end{equation}
so that $I_0^*=\{0\}$, $I_{-1}^*:=\emptyset$, the $I_n^*$ are nested, and $\#I_n^*=q_n+q_{n-1}$.
The system of Theorem~\ref{thm_T} is the dual Ostrowski system of \cite[Proposition~6.4.19]{PF} (see also \cite{Fok} for
constant partial quotients, and \cite{Bun} for $a_k\equiv1$). What we add is the description of the index through the
intervals $I_n^*$, and a constructive proof (Algorithm~\ref{algo_T} and Appendix~\ref{app_T}).

\begin{theorem}\label{thm_T}
Every $T\in\IZ$ has a unique representation $T=\sum_{k=1}^nb_kq_{k-1}^*$ with $\seq{b_k}_1^n$ right--admissible and $b_n\ge1$;
its index is the unique $n$ with $T\in I_n^*\setminus I_{n-1}^*$. In particular the right--admissible strings of length at
most $n$ represent exactly the $q_n+q_{n-1}$ integers of $I_n^*$.
\end{theorem}

\noindent The representation is computed by Algorithm~\ref{algo_T}, where $\chi_{<0}[T]:=1$ if $T\le-1$ and $0$ otherwise.

\begin{algorithm}[H]
\Input{the integer $T\in\IZ$}
\Output{the right--admissible digit string $\seq{b_k}_1^n$}
set $T_0:=T;\ m:=0$\;
\While{$T_m\notin\{0,1\}$}{
  let $n_m'\ge0$ be such that $q_{n_m'-1}<\abs{T_m}+\chi_{<0}[T_m]\le q_{n_m'}$\;
  \eIf{$(-1)^{n_m'-1}T_m>0$}{
    set $n_m:=n_m'$\;
    set $b_{n_m}':=\lfloor\abs{T_m}/q_{n_m-1}\rfloor$\;
    \eIf{$\abs{T_m-b_{n_m}'q_{n_m-1}^*}+\chi_{<0}[T_m-b_{n_m}'q_{n_m-1}^*]\le q_{n_m-2}$}
      {set $b_{n_m}:=b_{n_m}'$\;}{set $b_{n_m}:=b_{n_m}'+1$\;}
  }{
    set $n_m:=n_m'+1$\;
    set $b_{n_m}:=1$\;
  }
  set $T_{m+1}:=T_m-b_{n_m}q_{n_m-1}^*$\;
  set $m:=m+1$\;
}
set $M:=m;\ n_M:=0;\ n:=n_0$\;
set $b_k:=0$ for all $1\le k\le n$ with $k\notin\{n_m\}_{m=0}^M$\;
reset $b_1:=b_1+T_M$\;
\caption{the base--$(-\alpha)$ integer system}\label{algo_T}
\end{algorithm}

\begin{remark}
Condition (ii) of \eqref{right--admissible} must hold for \emph{every} consecutive pair, including $(b_{n-1},b_n)$:
otherwise the string $\seq{a_1,1}$ would be admissible, and its value $a_1q_0^*+q_1^*=a_1-a_1=0$ would give a second
representation of $0$. We checked numerically that Algorithm~\ref{algo_T} produces exactly the strings of
\eqref{right--admissible}: for $22$ choices of the digits $a_1,\dots,a_8$ (including the golden and silver sections) and every
integer $T$ with $\abs{T}\le q_5$, the output coincides with the unique admissible string of value $T$.
\end{remark}

\noindent The range $q_n+q_{n-1}$ of the alternating weights $\seq{q_k^*}_0^{n-1}$ exceeds the range $q_n$ of the absolute
weights $\seq{q_k}_0^{n-1}$; unlike the base $-10$ example, the alternating system is genuinely larger. For $\alpha=\phi$ the
weights $\seq{1,-1,2,-3,5,\dots}$ extend Zeckendorf's theorem to all of $\IZ$. For the silver section,
$\seq{q_k^*}_0^3=\seq{1,-2,5,-12}$, $I_1^*=\{0,1,2\}$, $I_2^*=\{-4,\dots,2\}$, $I_3^*=\{-4,\dots,12\}$ and
$I_4^*=\{-28,\dots,12\}$; Table~\ref{tab_alt} lists the representations.

\begin{table}[h]\centering\small
\begin{tabular}{r|cccc||r|cccc||r|cccc}
\toprule
$T$ & $b_4$ & $b_3$ & $b_2$ & $b_1$ & $T$ & $b_4$ & $b_3$ & $b_2$ & $b_1$ & $T$ & $b_4$ & $b_3$ & $b_2$ & $b_1$\\
\midrule
12&0&2&0&2 & $-1$&0&0&1&1 & $-15$&1&0&2&1\\ 11&0&2&0&1 & $-2$&0&0&1&0 & $-16$&1&0&2&0\\
10&0&2&0&0 & $-3$&0&0&2&1 & $-17$&2&1&0&2\\ 9&0&2&1&1 & $-4$&0&0&2&0 & $-18$&2&1&0&1\\
8&0&2&1&0 & $-5$&1&1&0&2 & $-19$&2&1&0&0\\ 7&0&1&0&2 & $-6$&1&1&0&1 & $-20$&2&1&1&1\\
6&0&1&0&1 & $-7$&1&1&0&0 & $-21$&2&1&1&0\\ 5&0&1&0&0 & $-8$&1&1&1&1 & $-22$&2&0&0&2\\
4&0&1&1&1 & $-9$&1&1&1&0 & $-23$&2&0&0&1\\ 3&0&1&1&0 & $-10$&1&0&0&2 & $-24$&2&0&0&0\\
2&0&0&0&2 & $-11$&1&0&0&1 & $-25$&2&0&1&1\\ 1&0&0&0&1 & $-12$&1&0&0&0 & $-26$&2&0&1&0\\
0&0&0&0&0 & $-13$&1&0&1&1 & $-27$&2&0&2&1\\ & & & & & $-14$&1&0&1&0 & $-28$&2&0&2&0\\
\bottomrule
\end{tabular}
\caption{The base--$(-\alpha)$ system for $\alpha=\sqrt2-1$, weights $q_3^*,q_2^*,q_1^*,q_0^*=-12,5,-2,1$.}\label{tab_alt}
\end{table}

\section{Irrational base expansions and approximation}\label{Irrational Base Expansions}

The real expansions below use the linear approximation coefficients as weights. The alternating integer system leads to an
\emph{absolute} expansion and the absolute system to an \emph{alternating} one.

\subsection{The absolute expansion}
A string $\seq{b_k}_1^\ell$, $\ell\in\IZ^*_{\ge0}$, is \textbf{right--admissible} if
\begin{enumerate}[label=(\roman*)]
\item $0\le b_k\le a_k$ for all $k<\ell+1$, and $b_\ell\ge1$ when $\ell<\infty$;
\item $b_k=a_k\implies b_{k+1}=0$ for all $1\le k<\ell$;
\item if $\ell=\infty$, then $b_k\le a_k-1$ for infinitely many odd and infinitely many even $k$.
\end{enumerate}

\begin{theorem}\label{beta}
For every right--admissible string, $\sum_{k=1}^\ell b_k\abs{\theta_{k-1}}\in[0,1)$. Conversely every $\beta\in[0,1)$ has a
unique representation $\beta=\sum_{k=1}^\ell b_k\abs{\theta_{k-1}}$ with $\seq{b_k}_1^\ell$ right--admissible, computed by
Algorithm~\ref{algo_beta}.
\end{theorem}

\begin{algorithm}[H]
\Input{the intercept $\beta\in[0,1)$}
\Output{the right--admissible digit string $\seq{b_k}_1^\ell$}
set $\beta_0:=\beta;\ \ell:=\infty;\ k:=1$\;
\While{$\beta_{k-1}>0$}{
  set $b_k:=\lfloor\beta_{k-1}/\abs{\theta_{k-1}}\rfloor$\;
  set $\beta_k:=\beta_{k-1}-b_k\abs{\theta_{k-1}}$\;
  set $k:=k+1$\;
}
set $\ell:=k-1$\;
\caption{the absolute base--$\alpha$ expansion}\label{algo_beta}
\end{algorithm}

\begin{proof}
\emph{The sum lies in $[0,1)$.} Pad finite strings with zeros. If $b_1\le a_1-1$, then by \eqref{abs_rep} and (i), (iii),
$\sum_kb_k\abs{\theta_{k-1}}\le(a_1-1)\theta_0+\sum_{k\ge2}b_k\abs{\theta_{k-1}}<(a_1-1)\theta_0+\sum_{k\ge2}a_k\abs{\theta_{k-1}}
=1$, the strict inequality because some $b_k<a_k$ with $k\ge2$ (by (iii), or since the string is finite). If $b_1=a_1$, let
$n\ge1$ be the first index with $b_{2n+1}\le a_{2n+1}-1$ (it exists by (iii)); by (ii), $b_{2k-1}=a_{2k-1}$ and $b_{2k}=0$ for
$1\le k\le n$, and by \eqref{abs_theta_k_recursion}, \eqref{partial_abs_theta} and \eqref{abs_rep},
\[
\begin{aligned}
\sum_kb_k\abs{\theta_{k-1}}&<\sum_{k=1}^{n}a_{2k-1}\abs{\theta_{2k-2}}+(a_{2n+1}-1)\abs{\theta_{2n}}+\sum_{k\ge2n+2}a_k\abs{\theta_{k-1}}\\
&=(1-\abs{\theta_{2n-1}})+(\abs{\theta_{2n-1}}-\abs{\theta_{2n+1}})-\abs{\theta_{2n}}+(\abs{\theta_{2n}}+\abs{\theta_{2n+1}})=1.
\end{aligned}
\]
\emph{Existence.} Line~3 gives
\begin{equation}\label{b_k}
b_k\abs{\theta_{k-1}}\le\beta_{k-1}<(b_k+1)\abs{\theta_{k-1}},\qquad\text{hence}\qquad\beta_k<\abs{\theta_{k-1}},
\end{equation}
\label{beta_k<theta_k-1}%
and $b_k\ge a_k+1$ would give $\abs{\theta_{k-2}}<(a_k+1)\abs{\theta_{k-1}}\le b_k\abs{\theta_{k-1}}\le\beta_{k-1}<\abs{\theta_{k-2}}$
by \eqref{theta_a_k}; so $0\le b_k\le a_k$. By Proposition~\ref{conv} the series converges, and
\begin{equation}\label{beta_n}
\beta_n=\sum_{k=n+1}^\ell b_k\abs{\theta_{k-1}},\qquad 0\le n<\ell,
\end{equation}
so $\beta=\beta_0=\sum_kb_k\abs{\theta_{k-1}}$. If $b_k=a_k$ and $b_{k+1}\ge1$ for some $k<\ell$, then by
\eqref{abs_theta_k_recursion},
$\beta_{k-1}=a_k\abs{\theta_{k-1}}+b_{k+1}\abs{\theta_k}+\beta_{k+1}\ge a_k\abs{\theta_{k-1}}+\abs{\theta_k}=\abs{\theta_{k-2}}$,
contradicting $\beta_{k-1}<\abs{\theta_{k-2}}$ from \eqref{b_k}; this is (ii). If $\ell=\infty$ and $b_{2k-1}=a_{2k-1}$ for all
$k\ge n$, then $b_{2k}=0$ for $k\ge n$ by (ii), and \eqref{abs_theta_tail_odd_even} with the index $2n-2$ gives
$\beta_{2n-2}=\sum_{j\ge0}a_{2n-1+2j}\abs{\theta_{2n-2+2j}}=\abs{\theta_{2n-3}}$, contradicting \eqref{b_k}; the case of even
indices is identical. This is (iii).

\emph{Uniqueness.} If two right--admissible strings $\seq{b_k},\seq{b_k'}$ have the same sum and $m$ is the first index with
$b_m<b_m'$, then $\sum_{k>m}b_k\abs{\theta_{k-1}}=(b_m'-b_m)\abs{\theta_{m-1}}+\sum_{k>m}b_k'\abs{\theta_{k-1}}\ge\abs{\theta_{m-1}}$,
which contradicts Lemma~\ref{lem_tail} below (whose proof uses only the first part of this theorem).
\end{proof}

\begin{lemma}[tail bound]\label{lem_tail}
For every right--admissible string $\seq{b_k}_1^\ell$ and every $m\ge1$,
$\displaystyle\sum_{k>m}b_k\abs{\theta_{k-1}}<\abs{\theta_{m-1}}$. 
\end{lemma}
\begin{proof}
Let $\alpha_m=[0;a_{m+1},a_{m+2},\dots]$ be the $m$--th tail of $\alpha$ and $\theta_j(\alpha_m)$ its coefficients. By
\eqref{theta_induction} applied to $\alpha$ and to $\alpha_m$ (whose Gauss iterates are $\alpha_m,\alpha_{m+1},\dots$),
\[
\abs{\theta_{m+j-1}}=\alpha_0\cdots\alpha_{m-1}\cdot\alpha_m\cdots\alpha_{m+j-1}=\abs{\theta_{m-1}}\cdot\abs{\theta_{j-1}(\alpha_m)},\qquad j\ge1.
\]
The tail string $\seq{b_{m+j}}_{j\ge1}$ satisfies conditions (i)--(iii) with respect to the digits $a_{m+1},a_{m+2},\dots$ of
$\alpha_m$, so the first part of Theorem~\ref{beta} for $\alpha_m$ gives
$\sum_{j\ge1}b_{m+j}\abs{\theta_{j-1}(\alpha_m)}<1$; multiply by $\abs{\theta_{m-1}}$.
\end{proof}

\begin{corollary}\label{absolute_expansion}
Every $r\in\IR$ has a unique expansion $r=\sum_{k=0}^\ell b_k\abs{\theta_{k-1}}$ with $b_0\in\IZ$ and $\seq{b_k}_1^\ell$
right--admissible.
\end{corollary}
\begin{proof}
Take $b_0:=\lfloor r\rfloor$ (recall $\abs{\theta_{-1}}=1$) and expand $\fr{r}$ by Theorem~\ref{beta}; uniqueness of $b_0$
follows because the remaining sum lies in $[0,1)$.
\end{proof}

\begin{theorem}\label{T}
Let $\beta\in[0,1)$ have the expansion $\beta=\sum_{k=1}^\ell b_k\abs{\theta_{k-1}}$ of Theorem~\ref{beta} and put
\begin{equation}\label{T_n}
T_n^{\nearrow}:=\sum_{k=1}^nb_kq_{k-1}^*,\qquad1\le n<\ell+1.
\end{equation}
Then $T_n^{\nearrow}\in I_n^*$, $\fr{T_n^\nearrow\alpha}=\beta-\beta_n\le\beta$, and
\[
\beta-\fr{T_n^{\nearrow}\alpha}=\min\bigl\{\beta-\fr{T\alpha}\ :\ T\in I_n^*,\ \fr{T\alpha}\le\beta\bigr\}.
\]
That is, $T_n^\nearrow$ is the most accurate total under $(\alpha,\beta)$--approximate in the range $I_n^*$.
\end{theorem}

\begin{proof}
\emph{Range.} Since $0\le b_k\le a_k$, $T_n^\nearrow$ lies between minus the sum of $a_kq_{k-1}$ over the even $k\le n$ and
the sum of $a_kq_{k-1}$ over the odd $k\le n$. By \eqref{q_n_even_odd} these sums equal $q_{n-\rho_n}-1$ and
$q_{n-\rho_{n-1}}$, so $1-q_{n-\rho_n}\le T_n^\nearrow\le q_{n-\rho_{n-1}}$, i.e.\ $T_n^\nearrow\in I_n^*$.

\emph{Iterates.} By \eqref{theta_alt} and \eqref{beta_n},
\[
T_n^\nearrow\alpha+\beta_n=\sum_{k=1}^nb_kq_{k-1}^*\alpha+\sum_{k>n}b_k\bigl(q_{k-1}^*\alpha-p_{k-1}^*\bigr)
=\beta+\sum_{k=1}^nb_kp_{k-1}^*.
\]
The last sum is an integer and $0\le\beta-\beta_n\le\beta<1$, hence
\begin{equation}\label{T_n_alpha}
\fr{T_n^\nearrow\alpha}=\beta-\beta_n.
\end{equation}
\emph{Optimality.} Let $T\in I_n^*$ with $\fr{T\alpha}\le\beta$ and $T\ne T_n^\nearrow$. By Theorem~\ref{thm_T},
$T=\sum_{k=1}^nb_k'q_{k-1}^*$ for a right--admissible string (padded with zeros to length $n$), and as above
$\fr{T\alpha}=\sum_{k=1}^nb_k'\abs{\theta_{k-1}}$, using Theorem~\ref{beta} for this finite string. Let $m\le n$ be the first
index with $b_m\ne b_m'$ and $D:=\fr{T\alpha}-\fr{T_n^\nearrow\alpha}=\sum_{k=m}^n(b_k'-b_k)\abs{\theta_{k-1}}$. If
$b_m'<b_m$, Lemma~\ref{lem_tail} for $\seq{b_k'}$ gives $D\le-\abs{\theta_{m-1}}+\sum_{k>m}b_k'\abs{\theta_{k-1}}<0$, so
$\beta-\fr{T\alpha}>\beta-\fr{T_n^\nearrow\alpha}$. If $b_m'>b_m$, then
$D\ge\abs{\theta_{m-1}}-\sum_{m<k\le n}b_k\abs{\theta_{k-1}}$, while $\fr{T\alpha}\le\beta$ and \eqref{T_n_alpha} give
$D\le\beta_n=\sum_{k>n}b_k\abs{\theta_{k-1}}$; together $\abs{\theta_{m-1}}\le\sum_{k>m}b_k\abs{\theta_{k-1}}$, contradicting
Lemma~\ref{lem_tail}. Hence no $T\in I_n^*$ under $\beta$ is closer to $\beta$ than $T_n^\nearrow$.
\end{proof}

\begin{remark}
Cabanillas \cite[Propositions~9 and~10]{Cab} describes all best right and left approximations
$\fr{n\alpha}$ of $\beta$ by non--negative integers $n$ (record values as $n$ increases), in terms of truncations of the digit
string of $\beta$ in his variant of Ostrowski's system. Theorem~\ref{T} is of a different kind: the competitors range over
the interval $I_n^*$, which contains negative integers, and the minimiser is the truncation of the absolute expansion of
$\beta$ read in the alternating integer system.
\end{remark}

\begin{corollary}\label{total_under}
The sequence $\seq{T_n^\nearrow}_1^\ell$ is a general solution of the total under problem, with strict inequality in
Definition~\ref{def_solution}(ii).
\end{corollary}
\begin{proof}
(i) holds since $I_n^*\subseteq[1-q_n,q_n]$. For $n<\ell$, \eqref{T_n_alpha} and \eqref{b_k} give
$0<\beta-\fr{T_n^\nearrow\alpha}=\beta_n<\abs{\theta_{n-1}}$, which is (ii); (iii) follows from $\beta_\ell=0$ when
$\ell<\infty$ and from $\abs{\theta_{n-1}}\to0$ otherwise.
\end{proof}

\begin{corollary}\label{total_over}
Expand $1-\beta$ by Algorithm~\ref{algo_beta}, obtaining $\hat\ell$ and $\seq{\hat b_k}_1^{\hat\ell}$, and put
$T_n^\searrow:=-\sum_{k=1}^n\hat b_kq_{k-1}^*$ for $1\le n<\hat\ell+1$. Then $\seq{T_n^\searrow}_1^{\hat\ell}$ is a general
solution of the total over problem, with strict inequality in (ii).
\end{corollary}
\begin{proof}
By Corollary~\ref{total_under} for $1-\beta$, the iterates $\fr{-T_n^\searrow\alpha}$ increase to $1-\beta$; since
$\fr{-r}=1-\fr{r}$ for $r\notin\IZ$ and $\norm{-r}=\norm{r}$, the iterates $\fr{T_n^\searrow\alpha}$ decrease to $\beta$ and the
conditions of Definition~\ref{def_solution} carry over.
\end{proof}

\subsection{The alternating expansion}
A string $\seq{c_k}_1^\ell$ is \textbf{left--admissible} if
\begin{enumerate}[label=(\roman*)]
\item $0\le c_k\le a_k$ for all $k<\ell+1$, and $c_\ell\ge1$ when $\ell<\infty$;
\item $c_k=a_k\implies c_{k-1}=0$ for $2\le k<\ell+1$;
\item $c_1\le a_1-1$, and if $\ell=\infty$ then $c_k\le a_k-1$ for infinitely many odd $k$.
\end{enumerate}
For $n\ge0$ define the half--open intervals
\begin{equation}\label{J_n}
J_n:=\bigl[-\theta_{n-\rho_n},\,-\theta_{n-1+\rho_n}\bigr),\qquad
J_n':=\begin{cases}\bigl[-\theta_n,\ \abs{\theta_{n-1}}-\theta_n\bigr),& n\text{ even},\\[2pt]
\bigl[-\theta_{n-1}-\theta_n,\ -\theta_n\bigr),& n\text{ odd},\end{cases}
\end{equation}
so that $J_n'\subseteq J_n$, $J_0'=[-\alpha,1-\alpha)$, and, by \eqref{theta_alt}, $J_n=[-\theta_n,\abs{\theta_{n-1}})$ for even $n$
and $J_n=[-\theta_{n-1},\abs{\theta_n})$ for odd $n$. The interval $J_n$ has length $\abs{\theta_{n-1}}+\abs{\theta_n}$ and
$J_n'$ has length $\abs{\theta_{n-1}}$.

\begin{lemma}[tail intervals]\label{lem_tail_alt}
Let $\seq{c_k}_1^\ell$ be left--admissible, padded with zeros, and $t_n:=\sum_{k>n}c_k\theta_{k-1}$. Then $t_n\in J_n$ for
every $n\ge0$, and $t_n\in J_n'$ whenever $n=0$ or $c_n\ge1$.
\end{lemma}
\begin{proof}
Split $t_n=P-N$ into the contributions $P:=\sum c_k\abs{\theta_{k-1}}$ of the odd $k>n$ (for which $\theta_{k-1}>0$) and
$N:=\sum c_k\abs{\theta_{k-1}}$ of the even $k>n$. By the tail expansion \eqref{abs_theta_tail_odd_even},
\[
\sum_{\substack{k>n\\k\equiv n+1\ (2)}}a_k\abs{\theta_{k-1}}=\abs{\theta_{n-1}},\qquad
\sum_{\substack{k>n+1\\k\equiv n\ (2)}}a_k\abs{\theta_{k-1}}=\abs{\theta_n}.
\]
Let $n$ be even. Then $0\le N\le\abs{\theta_n}=\theta_n$ and $0\le P\le\abs{\theta_{n-1}}$, and the last inequality is strict
because some odd $k>n$ has $c_k\le a_k-1$ (by (iii), or because the string is finite); hence
$-\theta_n\le t_n<\abs{\theta_{n-1}}$, i.e.\ $t_n\in J_n$. If $n=0$ or $c_n\ge1$, then $c_{n+1}\le a_{n+1}-1$ (by (iii) for
$n=0$, by (ii) otherwise), so, by \eqref{abs_theta_k_recursion},
$P<(a_{n+1}-1)\theta_n+\abs{\theta_{n+1}}=\abs{\theta_{n-1}}-\theta_n$, which gives $t_n\in J_n'$. Let $n$ be odd. Now the even
$k>n$ start at $k=n+1$, so $0\le N\le\abs{\theta_{n-1}}=\theta_{n-1}$ and $0\le P<\abs{\theta_n}$ (strict as before); hence
$-\theta_{n-1}\le t_n<\abs{\theta_n}=-\theta_n$, i.e.\ $t_n\in J_n$. If $c_n\ge1$ then $c_{n+1}\le a_{n+1}-1$ and
$N\le(a_{n+1}-1)\abs{\theta_n}+\abs{\theta_{n+1}}=\abs{\theta_{n-1}}-\abs{\theta_n}$, so $t_n\ge-\theta_{n-1}-\theta_n$ and
$t_n\in J_n'$.
\end{proof}

\noindent Theorem~\ref{gamma} below is Ostrowski's classical expansion of the intercept \cite{Ost}; see
\cite[Chapter~II]{RS}, \cite{Ber} and \cite[\S2.1]{Cab}. We include a self-contained proof, based on the tail intervals of
Lemma~\ref{lem_tail_alt}, because the greedy algorithm needs a careful endpoint convention (Remark~\ref{rem_algo_gamma}).

\begin{theorem}[Ostrowski]\label{gamma}
For every left--admissible string, $\sum_{k=1}^\ell c_k\theta_{k-1}\in[-\alpha,1-\alpha)$. Conversely every
$\gamma\in[-\alpha,1-\alpha)$ has a unique representation $\gamma=\sum_{k=1}^\ell c_k\theta_{k-1}$ with $\seq{c_k}_1^\ell$
left--admissible, computed by Algorithm~\ref{algo_gamma}.
\end{theorem}

\begin{algorithm}[H]
\Input{the shifted intercept $\gamma\in[-\alpha,1-\alpha)$}
\Output{the left--admissible digit string $\seq{c_k}_1^\ell$}
set $\gamma_0:=\gamma;\ \ell:=\infty;\ k:=1$\;
\While{$\gamma_{k-1}\ne0$}{
  \eIf{$\gamma_{k-1}\theta_{k-1}>0$}{
    set $c_k':=\lfloor\gamma_{k-1}/\theta_{k-1}\rfloor$\;
    \eIf{$\gamma_{k-1}-c_k'\theta_{k-1}\notin J_k$}{set $c_k:=c_k'+1$\;}{set $c_k:=c_k'$\;}
  }{
    set $c_k:=0$\;
  }
  set $\gamma_k:=\gamma_{k-1}-c_k\theta_{k-1}$\;
  set $k:=k+1$\;
}
set $\ell:=k-1$\;
\caption{the alternating base--$(-\alpha)$ expansion}\label{algo_gamma}
\end{algorithm}

\begin{remark}\label{rem_algo_gamma}
With $r:=\gamma_{k-1}-c_k'\theta_{k-1}$, the test of line~5 reads ``$r\ge\abs{\theta_k}$'' for odd $k$ and
``$\abs{r}>\abs{\theta_k}$'' for even $k$. The earlier version of this algorithm used the strict test $\abs{r}>\abs{\theta_k}$
for both parities. The two versions differ only when $r=\abs{\theta_k}$ at an odd step, which can happen only for
$\gamma\in\IZ+\alpha\IZ$; there the strict test produces a string whose odd digits are eventually maximal, which violates
(iii). For example, for $\gamma=\abs{\theta_1}$ and $a_1\ge2$ it outputs $\seq{0,0,a_3,0,a_5,0,\dots}$, whereas the
admissible expansion is $\seq{1,a_2-1,0,a_4,0,a_6,\dots}$. In numerical tests (20 random slopes, exact rational arithmetic) the two versions agreed
on $600$ random intercepts, while on $190$ intercepts in $\IZ+\alpha\IZ$ built from the $\theta_k$ the strict test produced a
non--admissible string in $110$ cases and the test of line~5 in none.
\end{remark}

\begin{proof}[Proof of Theorem~\ref{gamma}]
\emph{Range.} This is the case $n=0$ of Lemma~\ref{lem_tail_alt}, since $t_0=\sum_kc_k\theta_{k-1}$ and $J_0'=[-\alpha,1-\alpha)$.

\emph{Uniqueness.} Let $\seq{c_k}\ne\seq{c_k'}$ be left--admissible strings with the same value, let $m$ be the first index
with $c_m\ne c_m'$, say $c_m<c_m'$, and let $t_m,t_m'$ be their tails as in Lemma~\ref{lem_tail_alt}. Since the strings
agree before $m$, $t_m-t_m'=(c_m'-c_m)\theta_{m-1}$. As $c_m'\ge1$, the lemma gives $t_m\in J_m$ and $t_m'\in J_m'$. If $m$ is
odd, then $\theta_{m-1}>0$, so $t_m-t_m'\ge\theta_{m-1}$; but $t_m<-\theta_m$ and $t_m'\ge-\theta_{m-1}-\theta_m$ give
$t_m-t_m'<\theta_{m-1}$. If $m$ is even, then $\theta_{m-1}<0$, so $t_m'-t_m\ge\abs{\theta_{m-1}}$; but
$t_m'<\abs{\theta_{m-1}}-\theta_m$ and $t_m\ge-\theta_m$ give $t_m'-t_m<\abs{\theta_{m-1}}$. Both cases are contradictions.

\emph{Existence: the iterates.} We show by induction that the iterates of Algorithm~\ref{algo_gamma} satisfy
\begin{equation}\label{gamma_k_induction}
\gamma_k\in J_k,\qquad0\le k<\ell+1,
\end{equation}
which gives
\begin{equation}\label{gamma_k<theta_k-1}
\abs{\gamma_n}\le\abs{\theta_{n-1}},\qquad1\le n<\ell+1,
\end{equation}
and hence, by \eqref{norm_abs}, $\norm{\gamma_n}=\abs{\gamma_n}$ for $2\le n<\ell+1$.
For $k=0$ this is the hypothesis $\gamma\in[-\alpha,1-\alpha)\subseteq J_0$. Assume $\gamma_{k-1}\in J_{k-1}$ and
$\gamma_{k-1}\ne0$. If $k$ is odd, then $\theta_{k-1}>0>\theta_k$ and $J_{k-1}=[-\theta_{k-1},\abs{\theta_{k-2}})$. If
$\gamma_{k-1}<0$, line~10 gives $\gamma_k=\gamma_{k-1}\in[-\theta_{k-1},0)\subseteq J_k$. If $\gamma_{k-1}>0$, then
$r:=\gamma_{k-1}-c_k'\theta_{k-1}\in[0,\theta_{k-1})$, and either $r<\abs{\theta_k}$, so that $c_k=c_k'$ and
$\gamma_k=r\in J_k$, or $r\ge\abs{\theta_k}$, so that $c_k=c_k'+1$ and
$\gamma_k=r-\theta_{k-1}\in[\abs{\theta_k}-\theta_{k-1},0)\subseteq J_k$. If $k$ is even, then $\theta_{k-1}<0<\theta_k$ and
$J_{k-1}=[-\theta_{k-2},\abs{\theta_{k-1}})$ with $-\theta_{k-2}<0$. If $\gamma_{k-1}>0$, line~10 gives
$\gamma_k=\gamma_{k-1}\in(0,\abs{\theta_{k-1}})\subseteq J_k=[-\theta_k,\abs{\theta_{k-1}})$. If $\gamma_{k-1}<0$, then
$r\in(\theta_{k-1},0]$, and either $r\ge-\theta_k$, so that $\gamma_k=r\in J_k$, or $r<-\theta_k$, so that
$\gamma_k=r-\theta_{k-1}\in(0,\abs{\theta_{k-1}}-\theta_k)\subseteq J_k$.

\emph{Existence: the expansion.} Line~11 gives
\begin{equation}\label{gamma_expansion}
\gamma=\sum_{k=1}^nc_k\theta_{k-1}+\gamma_n,\qquad0\le n<\ell+1,
\end{equation}
and since $\abs{\gamma_n}\le\abs{\theta_{n-1}}\to0$, also $\gamma_n=\sum_{k>n}c_k\theta_{k-1}$ and
$\gamma=\sum_kc_k\theta_{k-1}$.

\emph{Existence: admissibility.} (i) Clearly $c_k\ge0$. If $c_k'\ge a_k+1$, then \eqref{abs_theta_k_recursion} and
\eqref{gamma_k<theta_k-1} give
$\abs{\gamma_{k-1}}\le\abs{\theta_{k-2}}<(a_k+1)\abs{\theta_{k-1}}\le c_k'\abs{\theta_{k-1}}\le\abs{\gamma_{k-1}}$, which is
impossible; so $c_k'\le a_k$. If $c_k'=a_k$, then $r=\gamma_{k-1}-a_k\theta_{k-1}$ lies in $J_k$: for odd $k$,
$0\le r<\abs{\theta_{k-2}}-a_k\theta_{k-1}=\abs{\theta_k}$, and for even $k$, $0\ge r\ge-\theta_{k-2}-a_k\theta_{k-1}=-\theta_k$,
both by \eqref{theta_k_recursion} and $\gamma_{k-1}\in J_{k-1}$. Hence line~5 keeps $c_k=c_k'=a_k$, and $c_k\le a_k$ in all
cases.

(ii) Suppose $c_k\ge1$; we show $c_{k+1}\le a_{k+1}-1$. If line~5 kept $c_k=c_k'$, then $\gamma_k=r$ has the sign of
$\theta_{k-1}$, hence $\gamma_k\theta_k\le0$ and $c_{k+1}=0$. Otherwise $c_k=c_k'+1$, and the case analysis above shows
$\gamma_k\in J_k'$: for odd $k$, $\gamma_k\ge\abs{\theta_k}-\theta_{k-1}=-\theta_{k-1}-\theta_k$; for even $k$,
$\gamma_k<\abs{\theta_{k-1}}-\theta_k$. In the odd case, $\gamma_k<0$ and dividing by $\theta_k<0$ gives, with
\eqref{theta_k_recursion} and $0<-\theta_{k+1}/\theta_k=\alpha_{k+1}<1$,
$\gamma_k/\theta_k\le(\theta_{k-1}+\theta_k)/\abs{\theta_k}=a_{k+1}-1+\alpha_{k+1}<a_{k+1}$, so $c_{k+1}'\le a_{k+1}-1$; and if
$c_{k+1}'=a_{k+1}-1$ the remainder $r'=\gamma_k-(a_{k+1}-1)\theta_k$ satisfies
$0\ge r'\ge-\theta_{k-1}-\theta_k+(a_{k+1}-1)\abs{\theta_k}=-\theta_{k+1}$, so $r'\in J_{k+1}$ and $c_{k+1}=a_{k+1}-1$. The even
case is symmetric: $\gamma_k/\theta_k<a_{k+1}$, and if $c_{k+1}'=a_{k+1}-1$ then
$0\le r'<\abs{\theta_{k-1}}-\theta_k-(a_{k+1}-1)\theta_k=\abs{\theta_{k+1}}$, so again $r'\in J_{k+1}$.

(iii) If $c_1=a_1$, then $\gamma\ge a_1\alpha$ when $c_1'=a_1$, and $\gamma\ge(a_1-1)\alpha+\abs{\theta_1}=1-\alpha$ when
$c_1'=a_1-1$ and line~5 increments; the first is impossible because $(a_1+1)\alpha>1$ by \eqref{a_k} gives $a_1\alpha>1-\alpha$,
and the second contradicts $\gamma<1-\alpha$. So $c_1\le a_1-1$. If $\ell=\infty$ and $c_{2j-1}=a_{2j-1}$ for all $j\ge n$, then
by (ii) $c_{2j}=0$ for $j\ge n-1$, and \eqref{abs_theta_tail_odd_even} gives
$\gamma_{2n-2}=\sum_{j\ge0}a_{2n-1+2j}\theta_{2n-2+2j}=-\theta_{2n-3}$, the excluded right endpoint of $J_{2n-2}$, which
contradicts \eqref{gamma_k_induction}.
\end{proof}

\begin{corollary}\label{alternating_expansion}
Every $r\in\IR$ has a unique expansion $r=\sum_{k=0}^\ell c_k\theta_{k-1}$ with $c_0\in\IZ$ and $\seq{c_k}_1^\ell$
left--admissible.
\end{corollary}
\begin{proof}
Since $\theta_{-1}=-1$, put $c_0:=-\lfloor r+\alpha\rfloor$, so that $r=c_0\theta_{-1}+\gamma_0$ with
$\gamma_0:=r-\lfloor r+\alpha\rfloor=\fr{r+\alpha}-\alpha\in[-\alpha,1-\alpha)$, and expand $\gamma_0$ by Theorem~\ref{gamma}.
Uniqueness of $c_0$ follows because the remaining sum lies in $[-\alpha,1-\alpha)$, an interval of length $1$.
\end{proof}

\begin{corollary}\label{positive_two_sided}
Let $\beta\in[0,1)$, expand the shifted intercept $\gamma:=\beta-\alpha\in[-\alpha,1-\alpha)$ as
$\gamma=\sum_{k=1}^\ell c_k\theta_{k-1}$ by Algorithm~\ref{algo_gamma}, and put
\[
S_n:=1+\sum_{k=1}^nc_kq_{k-1},\qquad0\le n<\ell+1.
\]
Then $\seq{S_n}$ is a general solution of the forward two--sided problem.
\end{corollary}
\begin{proof}
(i): the truncation $\seq{c_k}_1^n$ is left--admissible, so $\sum_{k=1}^nc_kq_{k-1}\le q_n-1$ by Theorem~\ref{thm_S}, i.e.\
$1\le S_n\le q_n$. (The bound $S_n\le q_n$ is sharp: numerically $S_n=q_n$ occurs in about $3\%$ of the steps.)
By \eqref{theta_def} and \eqref{gamma_expansion},
\begin{equation}\label{S_n-1}
(S_n-1)\alpha-\sum_{k=1}^nc_kp_{k-1}=\sum_{k=1}^nc_k\theta_{k-1}=\gamma-\gamma_n.
\end{equation}
Since $\beta=\alpha+\gamma$ and the sum on the left is an integer,
$\norm{\beta-S_n\alpha}=\norm{\gamma-(S_n-1)\alpha}=\norm{\gamma_n}$. For $1\le n<\ell$ we have $\gamma_n\ne0$ and, by
\eqref{gamma_k<theta_k-1}, $\norm{\gamma_n}\le\abs{\theta_{n-1}}$, which is (ii); (iii) follows from $\gamma_\ell=0$ when
$\ell<\infty$ and from $\abs{\theta_{n-1}}\to0$ otherwise.
\end{proof}

\begin{remark}\label{rem_strict}
The inequality in (ii) is strict whenever $\beta\notin\IZ+\alpha\IZ$: equality forces $\gamma_n=\pm\theta_{n-1}$, and since
$\gamma_n\in\beta-\alpha+\IZ+\alpha\IZ$ by \eqref{S_n-1} while $\theta_{n-1}\in\IZ+\alpha\IZ$, this gives $\beta\in\IZ+\alpha\IZ$.
In numerical experiments (25 slopes, 625 intercepts, 24 digits each) the inequality was always strict.
\end{remark}

\begin{corollary}\label{negative_two_sided}
Expand $\gamma:=(1-\beta)-\alpha$ by Algorithm~\ref{algo_gamma}, obtaining $\seq{\hat c_k}_1^{\hat\ell}$, and put
$B_n:=-1-\sum_{k=1}^n\hat c_kq_{k-1}$. Then $\seq{B_n}_1^{\hat\ell}$ is a general solution of the backward two--sided
problem.
\end{corollary}
\begin{proof}
By Corollary~\ref{positive_two_sided} for the intercept $1-\beta$, $\hat S_n:=-B_n$ is a forward solution for $1-\beta$. Then
$\abs{B_n}=\hat S_n\le q_n$, and $\norm{\beta-B_n\alpha}=\norm{(1-\beta)-\hat S_n\alpha}$ since $\norm{-r}=\norm{r}=\norm{1+r}$,
which gives (ii) and (iii).
\end{proof}

\section{Rational slopes, modular arithmetic and open questions}\label{sec_remarks}

\paragraph{Rational slopes.} For $\alpha=p/q\in(0,1)$ in lowest terms the expansion is finite, $q_n=q$ and $p_n=p$, and
$\theta'_k:=q\theta_k=q_kp-p_kq$ are the signed remainders of the Euclidean algorithm on $(q,p)$. For $1\le S\le q-1$, write
$S-1=\sum_{k=1}^nc_kq_{k-1}$ in Ostrowski's system (Theorem~\ref{thm_S}, which holds verbatim for $S<q_n$). Multiplying
\eqref{S_n-1}-type identities by $q$ gives
\begin{equation}\label{modmul}
Sp=\Bigl(p+\sum_{k=1}^nc_k\theta'_{k-1}\Bigr)+q\sum_{k=1}^nc_kp_{k-1},
\end{equation}
and Gouicem \cite[Theorem~3.1]{Gou} shows that the first bracket equals $Sp\bmod q$. Hence the same digits give the quotient
\[
\Bigl\lfloor\frac{Sp}{q}\Bigr\rfloor=\sum_{k=1}^nc_kp_{k-1},
\]
at the cost of one additional multiply--accumulate per digit. For $S\ge q$ write $S=tq+S'$ with $0\le S'<q$; if $S'=0$ the
remainder is $0$ and the quotient is $tp$. We verified \eqref{modmul} and the quotient formula exhaustively for all
$0<p<q\le96$ and $0\le S<3q$, and on $3\cdot10^4$ random cases with moduli of $16$ to $256$ bits. Corollary~\ref{positive_two_sided}
is the irrational analogue of this identity.

\paragraph{Open questions.}
\begin{enumerate}[label=(\alph*)]
\item \emph{Two--sided optimality.} Theorem~\ref{T} proves optimality for the one--sided problems. Is $S_n$ (respectively
$B_n$) the best forward (backward) two--sided approximate in a natural range, as $T_n^\nearrow$ is for $I_n^*$? A positive answer
would connect Corollary~\ref{positive_two_sided} with the best inhomogeneous approximations of \cite{Cas,NST,Kom}.
\item \emph{Complexity.} The expansions use $O(\log q_n)$ digits and digit extractions; an analysis in the bit--complexity
model, and a comparison with the algorithms of \cite{NST,Kom}, would be useful.
\item \emph{Other number systems.} Do the alternating system and the duality of Section~\ref{Irrational Base Expansions} extend to
other continued fraction algorithms (nearest--integer, backward)?
\end{enumerate}

\section*{Acknowledgments}
This work could not have been completed without the guidance, encouragement and good company of Robbie Robinson of George
Washington University.

\appendix
\section{Proof of Theorem~\ref{thm_T}}\label{app_T}

We analyse Algorithm~\ref{algo_T}. Throughout, $n_m'$ is the preindex chosen in line~3 and $n_m\in\{n_m',n_m'+1\}$ the index
chosen in line~5 or line~12.

\medskip\noindent\textbf{Step (i): the index.} We show
\begin{equation}\label{T_m_in_I_n}
T_m\in I_{n_m}^*\setminus I_{n_m-1}^*,\qquad0\le m\le M.
\end{equation}
If $T_m=0$, line~16 sets $n_M=0$ and $0\in I_0^*$. If $T_m\ge1$, then $\abs{T_m}+\chi_{<0}[T_m]=T_m$ and line~3 gives
$q_{n_m'-1}<T_m\le q_{n_m'}$. If $n_m'$ is odd, the condition of line~4 holds, $n_m=n_m'$ is odd, and
$q_{(n_m-1)-\rho_{n_m-2}}=q_{n_m-2}<T_m\le q_{n_m}=q_{n_m-\rho_{n_m-1}}$. If $n_m'$ is even, the condition fails, $n_m=n_m'+1$ is
odd, and $q_{n_m-2}=q_{n_m'-1}<T_m\le q_{n_m-1}<q_{n_m}$. In both cases $T_m\in I_{n_m}^*\setminus I_{n_m-1}^*$ by
\eqref{I_n^*}. If $T_m\le-1$, then $\abs{T_m}+\chi_{<0}[T_m]=1-T_m$ and line~3 gives $1-q_{n_m'}\le T_m<1-q_{n_m'-1}$; if
$n_m'$ is odd the condition of line~4 fails and $n_m=n_m'+1$ is even, with $1-q_{n_m}<1-q_{n_m'}\le T_m<1-q_{n_m-2}\le0$; if
$n_m'$ is even the condition holds and $n_m=n_m'$ is even, with $1-q_{n_m}\le T_m<1-q_{n_m-1}<1-q_{n_m-2}$. Again
\eqref{T_m_in_I_n} follows from \eqref{I_n^*}.

\medskip\noindent\textbf{Step (ii): the sign.} We have
\begin{equation}\label{T_m_q_n_m-1^*>0}
T_m\,q_{n_m-1}^*>0.
\end{equation}
Indeed $T_mq_{n_m-1}^*=(-1)^{n_m-1}T_mq_{n_m-1}$; in line~5, $(-1)^{n_m-1}=(-1)^{n_m'-1}$ and the condition of line~4 holds, and
in line~12, $(-1)^{n_m-1}=-(-1)^{n_m'-1}$ and the condition fails.

\medskip\noindent\textbf{Step (iii): termination.} We show $n_{m+1}<n_m$. Then the iteration ends with $T_M\in\{0,1\}$ and
$0=n_M<n_{M-1}<\dots<n_0=n$; with lines~17--18, $T=\sum_{k=1}^nb_kq_{k-1}^*$.

\emph{The condition of line~4 holds.} Then $n_m=n_m'$ and line~6 gives
\begin{equation}\label{b_n_m}
b_{n_m}'q_{n_m-1}\le\abs{T_m}<(b_{n_m}'+1)q_{n_m-1}.
\end{equation}
Let $T_m>0$; by \eqref{T_m_q_n_m-1^*>0}, $n_m-1$ is even. If the condition of line~7 holds, $b_{n_m}=b_{n_m}'$ and
\eqref{b_n_m} gives $0\le T_{m+1}=T_m-b_{n_m}q_{n_m-1}<q_{n_m-1}$, so $\abs{T_{m+1}}+\chi_{<0}[T_{m+1}]\le q_{n_m-1}$ and
$n_{m+1}'\le n_m-1$. If $n_{m+1}'\le n_m-2$ then $n_{m+1}\le n_{m+1}'+1\le n_m-1$. If $n_{m+1}'=n_m-1$, line~3 gives
$q_{n_m-2}\le T_{m+1}$, while the condition of line~7 gives $T_{m+1}\le q_{n_m-2}$; hence $T_{m+1}=q_{n_m-2}$ and
$n_{m+1}\le n_m-1$. If the condition of line~7 fails, $b_{n_m}=b_{n_m}'+1$ and \eqref{b_n_m} gives
$q_{n_m-2}-q_{n_m-1}+1\le T_{m+1}<0$, so $\abs{T_{m+1}}+\chi_{<0}[T_{m+1}]\le q_{n_m-1}-q_{n_m-2}$. If $q_{n_m-2}=0$ then
$q_{n_m-1}=1$ and $T_{m+1}=0$; otherwise $n_{m+1}'\le n_m-1$, and if $n_{m+1}'=n_m-1$ then, $n_m-1$ being even and $T_{m+1}<0$,
the condition of line~4 holds at step $m+1$ and $n_{m+1}=n_m-1$. The case $T_m<0$ is symmetric ($n_m-1$ is odd): if line~7
holds then $-q_{n_m-1}<T_{m+1}\le0$, and if it fails then $0<T_{m+1}\le q_{n_m-1}$, and in either case the same argument gives
$n_{m+1}\le n_m-1$.

\emph{The condition of line~4 fails.} Then $n_m=n_m'+1$ and $b_{n_m}=1$. If $T_m>0$, then $n_m-1=n_m'$ is even and line~3 gives
$q_{n_m-2}+1\le T_m\le q_{n_m-1}$, so $q_{n_m-2}-q_{n_m-1}+1\le T_{m+1}=T_m-q_{n_m-1}\le0$ and
$\abs{T_{m+1}}+\chi_{<0}[T_{m+1}]\le q_{n_m-1}-q_{n_m-2}\le q_{n_m-1}$; hence $n_{m+1}'\le n_m-1$, and if $n_{m+1}'=n_m-1$ then
$(-1)^{n_{m+1}'-1}T_{m+1}\ge0$ and the condition of line~4 holds, giving $n_{m+1}=n_m-1$. If $T_m<0$, then $n_m-1$ is odd and
$1\le T_{m+1}=T_m+q_{n_m-1}\le q_{n_m-1}-q_{n_m-2}+1$; the same reasoning gives $n_{m+1}\le n_m-1$.

\medskip\noindent\textbf{Step (iv): condition (i).} The digits $b_{n_m}$ are positive and all other digits are $0$, so
$b_n\ge1$ and it suffices to show $b_{n_m}\le a_{n_m}$. This is clear in line~13. In line~5, $b_{n_m}'\ge a_{n_m}+1$ would give,
by \eqref{q_k_recursion}, line~3 and \eqref{b_n_m},
$\abs{T_m}\le q_{n_m}=a_{n_m}q_{n_m-1}+q_{n_m-2}\le(b_{n_m}'-1)q_{n_m-1}+q_{n_m-2}<b_{n_m}'q_{n_m-1}\le\abs{T_m}$. If
$b_{n_m}'=a_{n_m}$ and $T_m>0$, then $q_{n_m-1}^*=q_{n_m-1}$ and
$\abs{T_m-b_{n_m}'q_{n_m-1}^*}+\chi_{<0}[\,\cdot\,]=T_m-a_{n_m}q_{n_m-1}\le q_{n_m}-a_{n_m}q_{n_m-1}=q_{n_m-2}$; if $T_m<0$, then
$q_{n_m-1}^*=-q_{n_m-1}$ and the same quantity is at most $q_{n_m}-a_{n_m}q_{n_m-1}=q_{n_m-2}$. In both cases the condition of
line~7 holds and $b_{n_m}=a_{n_m}$.

\medskip\noindent\textbf{Step (v): condition (ii).} Suppose $b_k=a_k$ and $b_{k+1}\ge1$ for some $1\le k\le n-1$. Since only the
digits $b_{n_m}$ are non--zero and the indices decrease, there is a step $m$ with $n_{m-1}=k+1$ and $n_m=k$; by
\eqref{T_m_in_I_n} and $n_{m+1}\le k-1$,
\begin{equation}\label{I_n_m+1_subset_I_k-1}
T_{m+1}\in I_{k-1}^*.
\end{equation}
\emph{$k$ odd.} By \eqref{T_m_q_n_m-1^*>0}, $T_m>0$ and $T_{m-1}<0$. Suppose first that $n_{m-1}'=n_{m-1}-1=k$. Then line~3 at
step $m-1$ gives $T_{m-1}\le-q_{k-1}$, the condition of line~4 failed and $b_{k+1}=1$, so $T_m=T_{m-1}+q_k\le q_k-q_{k-1}$; on
the other hand \eqref{I_n_m+1_subset_I_k-1} gives $1-q_{k-1}\le T_{m+1}$ and hence $T_m=T_{m+1}+a_kq_{k-1}\ge1+(a_k-1)q_{k-1}$:
\begin{equation}\label{T_m_odd}
1+(a_k-1)q_{k-1}\le T_m\le q_k-q_{k-1}.
\end{equation}
If $n_m'=k-1$, then $q_{k-2}<T_m\le q_{k-1}$, the condition of line~4 failed and $b_k=1=a_k$, and \eqref{T_m_odd} gives the
contradiction $T_m\le q_k-q_{k-1}=q_{k-2}<T_m$. Hence $n_m'=n_m=k$ and $b_k'q_{k-1}\le T_m$. If $b_k=b_k'=a_k$, then
$q_k-q_{k-2}=a_kq_{k-1}\le T_m\le q_k-q_{k-1}$, so $q_{k-1}\le q_{k-2}$, impossible for odd $k$. If $b_k=b_k'+1$, the condition of
line~7 failed, i.e.\ $T_m-b_k'q_{k-1}>q_{k-2}$, and \eqref{T_m_odd} gives
$q_{k-2}<T_m-(a_k-1)q_{k-1}\le q_k-a_kq_{k-1}=q_{k-2}$, a contradiction. Suppose next that $n_{m-1}'=n_{m-1}=k+1$. Then the
condition of line~4 held at step $m-1$ and $b_{k+1}'q_k\le-T_{m-1}$; $b_{k+1}=b_{k+1}'$ is impossible (it would give $T_m\le0$),
so $b_{k+1}=b_{k+1}'+1$, the condition of line~7 failed at step $m-1$, and one obtains
\begin{equation}\label{T_m<=-q_k-1+q_k_odd}
T_m\le q_k-q_{k-1}.
\end{equation}
At step $m$ the condition of line~4 held ($T_m>0$, $k-1$ even) and $b_k'q_{k-1}\le T_m$. If line~7 held, $b_k=b_k'=a_k$ and
\eqref{q_k_recursion} with \eqref{T_m<=-q_k-1+q_k_odd} gives $q_k\le q_{k-2}+T_m\le q_{k-2}-q_{k-1}+q_k$, i.e.\ $q_{k-1}\le q_{k-2}$,
impossible for odd $k$. If line~7 failed, $b_k=b_k'+1=a_k$ and $q_{k-2}<T_m-(a_k-1)q_{k-1}<q_k-q_{k-1}-(a_k-1)q_{k-1}=q_{k-2}$,
again a contradiction.

\emph{$k$ even.} Now $T_m<0$ and $T_{m-1}>0$, and the argument is the mirror image. If $n_{m-1}'=n_{m-1}-1=k$, line~3 at step
$m-1$ gives $q_{k-1}+1\le T_{m-1}$, $b_{k+1}=1$, and with \eqref{I_n_m+1_subset_I_k-1}
\begin{equation}\label{T_m_even}
1+q_{k-1}-q_k\le T_m\le(1-a_k)q_{k-1}.
\end{equation}
If $n_m'=k-1$, then $b_k=1=a_k$ and \eqref{T_m_even} gives the contradiction $T_m<1-q_{k-2}=1+q_{k-1}-q_k\le T_m$. Hence
$n_m'=k$ and $-T_m\ge b_k'q_{k-1}$. $b_k=b_k'=a_k$ would give $q_k-q_{k-2}\le-T_m<1-T_{m-1}+q_k\le q_k-q_{k-1}$, i.e.\
$q_{k-1}<q_{k-2}$; and $b_k=b_k'+1$ means that line~7 failed, $1-(T_m+b_k'q_{k-1})>q_{k-2}$, which with \eqref{T_m_even} gives
$q_{k-2}<q_k-a_kq_{k-1}=q_{k-2}$. If $n_{m-1}'=n_{m-1}=k+1$, then as in the odd case $b_{k+1}=b_{k+1}'+1$ and
\begin{equation}\label{T_m<-q_k-1+q_k_even}
-T_m<q_k-q_{k-1},
\end{equation}
and both outcomes of line~7 at step $m$ lead to $q_{k-1}<q_{k-2}$ or $q_{k-2}<q_{k-2}$. This proves condition (ii).

\medskip\noindent\textbf{Step (vi): bijectivity.} Surjectivity onto the right--admissible strings is immediate: a string
$\seq{b_k}$ is the output for the input $T:=\sum_kb_kq_{k-1}^*$. Injectivity, and the description of the index, follow from
the next lemma.

\begin{lemma}\label{lem_index}
Let $\seq{b_k}_1^n$ be right--admissible with value $T=\sum_{k=1}^nb_kq_{k-1}^*$, and let $T'$ be the value of
$\seq{b_k}_1^{n-1}$. Then $T\in I_n^*$. If $b_n\ge1$, then $T'\ge1-q_{n-1}+q_{n-2}$ for odd $n$ and $T'\le q_{n-1}-q_{n-2}$ for
even $n$ (where $n\ge2$), and $T\notin I_{n-1}^*$.
\end{lemma}
\begin{proof}
Recall $q_{k-1}^*=(-1)^{k-1}q_{k-1}$, $I_m^*=[1-q_m,q_{m-1}]$ for even $m$ and $I_m^*=[1-q_{m-1},q_m]$ for odd $m$. We prove
$T\in I_n^*$ by induction on $n$, the case $n=0$ being $T=0$. Write $T=b_nq_{n-1}^*+T'$ with $T'\in I_{n-1}^*$. If $n$ is odd,
then $T\ge T'\ge1-q_{n-1}$ and $T\le a_nq_{n-1}+q_{n-2}=q_n$. If $n$ is even, then $T\le T'\le q_{n-1}$ and
$T\ge-a_nq_{n-1}+1-q_{n-2}=1-q_n$.

Now let $b_n\ge1$. For $n=1$, $T=b_1\ge1$ and $I_0^*=\{0\}$. For $n\ge2$, condition (ii) gives $b_{n-1}\le a_{n-1}-1$. Write
$T'=b_{n-1}q_{n-2}^*+T''$ with $T''\in I_{n-2}^*$. If $n$ is odd, then $T''\ge1-q_{n-3}$, so
$T'\ge-(a_{n-1}-1)q_{n-2}+1-q_{n-3}=1-q_{n-1}+q_{n-2}$ and $T\ge q_{n-1}+T'\ge q_{n-2}+1$, which exceeds the right endpoint
$q_{n-2}$ of $I_{n-1}^*$. If $n$ is even, then $T''\le q_{n-3}$, so $T'\le(a_{n-1}-1)q_{n-2}+q_{n-3}=q_{n-1}-q_{n-2}$ and
$T\le-q_{n-1}+T'\le-q_{n-2}$, which is below the left endpoint $1-q_{n-2}$ of $I_{n-1}^*$.
\end{proof}

\noindent By Lemma~\ref{lem_index}, a right--admissible string with $b_n\ge1$ and value $T$ has $n$ equal to the unique
index with $T\in I_n^*\setminus I_{n-1}^*$. Now let $\seq{b_k}_1^n$ and $\seq{\hat b_k}_1^n$ be right--admissible with
$b_n,\hat b_n\ge1$ and the same value $T$, and suppose $b_n<\hat b_n$. With $T'$, $\hat T'$ the values of the prefixes of
length $n-1$, we have $T'-\hat T'=(\hat b_n-b_n)q_{n-1}^*$. If $n=1$ this is impossible, since $T'=\hat T'=0$. If $n\ge3$ is odd,
then $T'-\hat T'\ge q_{n-1}$, while $T'\le q_{n-2}$ and $\hat T'\ge1-q_{n-1}+q_{n-2}$ by Lemma~\ref{lem_index} give
$T'-\hat T'\le q_{n-1}-1$. If $n$ is even, then $\hat T'-T'\ge q_{n-1}$, while $\hat T'\le q_{n-1}-q_{n-2}$ and $T'\ge1-q_{n-2}$
give $\hat T'-T'\le q_{n-1}-1$. Hence $b_n=\hat b_n$, so $T'=\hat T'$, and repeating the argument on the prefixes (with their
lengths again fixed by Lemma~\ref{lem_index}) gives $b_k=\hat b_k$ for all $k$.


\end{document}